\documentclass[final,leqno]{siamltex}
\usepackage{amsmath}
\usepackage{amssymb}
\usepackage{graphicx}
\usepackage[notcite,notref]{showkeys}
\usepackage{tikz}
 \numberwithin{dummy}{section}

\newtheorem{algorithm}{Weak Galerkin Algorithm}

\newcommand{\bu}{{\bf u}}
\newcommand{\bq}{{\bf q}}

\newcommand{\be}{{\bf e}}
\newcommand{\bv}{{\bf v}}

\newcommand{\bpsi}{{\boldsymbol\psi}}

\def\Q{{\mathbb Q}}
\def\T{{\mathcal T}}
\def\E{{\mathcal E}}

\def\pT{{\partial T}}
\def\l{{\langle}}
\def\r{{\rangle}}

\def\T{{\mathcal T}}
\def\E{{\mathcal E}}

\def\bbf{{\bf f}}
\def\bg{{\bf g}}
\def\bn{{\bf n}}
\def\bq{{\bf q}}

\def\3bar{{|\hspace{-.02in}|\hspace{-.02in}|}}

 \def\an#1{\begin{align}#1\end{align}} 
 
\def\p#1{\begin{pmatrix}#1\end{pmatrix}}

\title{A stabilizer free WG Method for the Stokes Equations with order two superconvergence on polytopal mesh }
\author{Xiu Ye\thanks{Department of
Mathematics, University of Arkansas at Little Rock, Little Rock, AR
72204 (xxye@ualr.edu). This research was supported in part by
National Science Foundation Grant DMS-1620016.}
\and
Shangyou Zhang\thanks{Department of
Mathematical Sciences, University of Delaware, Newark, DE 19716 (szhang@udel.edu).}
}
\begin{document}
\maketitle

\begin{abstract}
A stabilizer free WG method is introduced for the Stokes equations with superconvergence on polytopal mesh in primary velocity-pressure formulation.
Convergence rates two order higher than the optimal-order for velocity of the WG approximation is proved in both an energy norm and the $L^2$ norm. Optimal order error estimate for pressure in the $L^2$ norm
is also established.
The numerical examples cover low and high order approximations, and 2D and 3D cases.
\end{abstract}


\pagestyle{myheadings}

\section{Introduction}
A stabilizing/penalty term  is often used in finite element methods
with discontinuous approximations  to enforce  connection of discontinuous functions
across element boundaries. Development of  stabilizer free discontinuous  finite element method is desirable since it  simplifies finite element
formulation and reduces programming complexity. The stabilizer free WG method and the stabilizer DG method on polytopal mesh were first introduced in \cite{sf-wg,cdg2} for second order elliptic problems.  The main idea in \cite{sf-wg,cdg2} is to raise the degree of polynomials used to compute weak gradient $\nabla_w$. In \cite{sf-wg, cdg2}, gradient is approximated by a polynomial of order $j=k+n-1$ where $n$ is the number of sides of  polygonal element.
This result has been improved in \cite{aw,aw1} by reducing the degree of polynomial $j$. Recently, new stabilizer free WG methods have been developed in \cite{sf-wg-part-II,sf-wg-part-III} for second order elliptic equations on polytopal mesh, which have superconvergence. Wachspress coordinates  are used to approximate $\nabla_w$ in \cite{liu, mu} for solving the Stokes equations on polytopal mesh. Wachspress coordinates are usually rational functions, instead of polynomials. The WG methods in \cite{liu, mu} are limited to the lowest order WG elements.

In this paper, we introduce a new stabilizer free WG method of any order to solve the Stokes problem: find unknown functions $\bu$ and $p$ such that
\begin{eqnarray}
-\Delta\bu+\nabla p&=& \bbf\quad
\mbox{in}\;\Omega,\label{moment}\\
\nabla\cdot\bu &=&0\quad \mbox{in}\;\Omega,\label{cont}\\
\bu&=&0\quad \mbox{on}\; \partial\Omega,\label{bc}
\end{eqnarray}
where $\Omega$ is a polygonal or polyhedral domain in
$\mathbb{R}^d\; (d=2,3)$.
Our new WG method has the following formulations without any stabilizers: seek $(\bu_h,p_h)\in V_h\times W_h$ satisfying the following for all $(\bv,w)\in V_h\times W_h$,
\begin{eqnarray}
(\nabla_w\bu_h,\nabla_w\bv)-(\nabla_w\cdot\bv, p_h)&=&({\bf f}, \bv),\label{w1}\\
(\nabla_w\cdot\bu_h, w)&=&0.\label{w2}
\end{eqnarray}
Here  $\nabla_w$ and $\nabla_w\cdot$ are weak gradient and weak divergence, respectively. In addition, we have proved that the WG approximations have the convergence rates two order higher than the optimal-order for velocity  in both an energy norm and the $L^2$ norm and the optimal convergence rate for pressure in the $L^2$ norm. Extensive numerical examples are tested for the new WG elements of different degrees $k$ in both two and three dimensional spaces.

\section{Preliminary}
Let ${\cal T}_h$ be a partition of the domain $\Omega$ consisting of
polygons in two dimension or polyhedra in three dimension satisfying
a set of conditions specified in \cite{wymix}. Denote by ${\cal E}_h$
the set of all edges or flat faces in ${\cal T}_h$, and let ${\cal
E}_h^0={\cal E}_h\backslash\partial\Omega$ be the set of all
interior edges or flat faces. For every element $T\in \T_h$, we
denote by $h_T$ its diameter and mesh size $h=\max_{T\in\T_h} h_T$
for ${\cal T}_h$. Let $P_k(T)$ consist all the polynomials on $T$ with degree no greater than $k$.

For $k\ge 0$ and given $\T_h$, define two finite element spaces  for velocity
\begin{eqnarray}
V_h &=&\left\{ \bv=\{\bv_0,\bv_b\}:\ \bv_0|_{T}\in [P_{k}(T)]^d,\;\bv_b|_e\in [P_{k+1}(e)]^d,e\subset\pT \right\},\label{vh}
\end{eqnarray}
and for pressure
\begin{equation}
W_h =\left\{w\in L_0^2(\Omega): \ w|_T\in P_{k+1}(T)\right\}.\label{wh}
\end{equation}
Let $V_h^0$ be a subspace of $V_h$ consisting of functions with vanishing boundary value.

The space $H(\operatorname{div};\Omega)$ is defined as 
\[
H(\operatorname{div}; \Omega)=\left\{ \bv\in [L^2(\Omega)]^d:\; \nabla\cdot\bv \in L^2(\Omega)\right\}.
\]
For any $T\in\T_h$, it can be divided in to a set of disjoint triangles $T_i$ with $T=\cup T_i$.  Then we define a space $\Lambda_h(T)$ for the approximation of weak gradient on each element $T$  as
\begin{eqnarray}
\Lambda_{k}(T)=\{\bpsi\in [H(\operatorname{div};T)]^d:&&\ \bpsi|_{T_i}\in [P_{k+1}(T_i)]^{d\times d},\;\;\nabla\cdot\bpsi\in [P_k(T)]^d,\label{lambda}\\
&&\bpsi\cdot\bn|_e\in [P_{k+1}(e)]^d,\;e\subset\pT\}.\nonumber
\end{eqnarray}
For a function $\bv\in V_h$, its weak gradient $\nabla_w\bv$ is a piecewise polynomial satisfying $\nabla_w\bv|_T \in \Lambda_k(T)$  and  the following equation,
\begin{equation}\label{wg}
(\nabla_w\bv,\  \tau)_T = -(\bv_0,\  \nabla\cdot \tau)_T+
\l\bv_b, \ \tau\cdot\bn \r_\pT\quad\forall\tau\in \Lambda_k(T).
\end{equation}
For a function $\bv\in V_h$, its weak divergence $\nabla_w\cdot\bv$ is a piecewise polynomial satisfying $\nabla_w\cdot\bv|_T \in P_{k+1}(T)$  and  the following equation,
\begin{equation}\label{wd}
(\nabla_w\cdot\bv,\  w)_T = -(\bv_0,\  \nabla w)_T+
\l\bv_b\cdot\bn, \ w \r_\pT\quad\forall w\in P_{k+1}(T).
\end{equation}

The proof of the following lemma can be found in \cite{sf-wg-part-III}.

\begin{lemma}
For $\tau\in [H(\operatorname{div};\Omega)]^d$, there exists a projection $\Pi_h$ with $\Pi_h\tau\in [H(\operatorname{div};\Omega)]^d$ satisfying $\Pi_h\tau|_T\in \Lambda_k(T)$ and the followings
\begin{eqnarray}
(\nabla\cdot\tau,\;\bq)_T&=&(\nabla\cdot\Pi_h\tau,\;\bq)_T\quad \forall\bq\in [P_k(T)]^d , \label{zhang1}\\
-(\nabla\cdot\tau, \;\bv_0)&=&(\Pi_h\tau, \;\nabla_w \bv)\quad\forall\bv=\{\bv_0,\bv_b\}\in V_h,\label{key}\\
\|\Pi_h\tau-\tau\|&\le& Ch^{k+2}|\tau|_{k+2}.\label{zhang2}
\end{eqnarray}
\end{lemma}

\smallskip


\section{ Finite Element Method and Its Well Posedness}

We start this section by introducing the following  WG finite element scheme without stabilizers.

\begin{algorithm}
A numerical approximation for (\ref{moment})-(\ref{bc}) is finding $(\bu_h,p_h)\in V_h^0\times W_h$  such that for all $(\bv,w)\in V_h^0\times  W_h$,
\begin{eqnarray}
(\nabla_w\bu_h,\ \nabla_w\bv)-(\nabla_w\cdot\bv,\;p_h)&=&(f,\;\bv),\label{wg1}\\
(\nabla_w\cdot\bu_h,\;w)&=&0.\label{wg2}
\end{eqnarray}
\end{algorithm}
Let $Q_0$ and $Q_b$ be the two element-wise defined $L^2$ projections onto $[P_k(T)]^d$ and $[P_{k+1}(e)]^d$ with $e\subset\partial T$ on $T$ respectively. Define $Q_h\bu=\{Q_0\bu,Q_b\bu\}\in V_h$ for the true solution $\bu$. Let $\Q_h$ be the element-wise defined $L^2$ projection onto $\Lambda_k(T)$ on each element $T$. Finally denote by $\mathcal{Q}_h$ the element-wise defined $L^2$ projection onto $P_{k+1}(T)$ on each element $T$.

\begin{lemma}
Let $\boldsymbol\phi\in [H_0^1(\Omega)]^d$, then on  $T\in\T_h$
\begin{eqnarray}
\nabla_w Q_h\boldsymbol\phi &=&\Q_h\nabla\boldsymbol\phi,\label{key1}\\
\nabla_w \cdot Q_h\boldsymbol\phi &=&\mathcal{Q}_h\nabla\cdot\boldsymbol\phi.\label{key2}
\end{eqnarray}
\end{lemma}
\begin{proof}
Using (\ref{wg}) and  integration by parts, we have that for
any $\tau\in \Lambda_k(T)$
\begin{eqnarray*}
(\nabla_w Q_h\boldsymbol\phi,\tau)_T &=& -(Q_0\boldsymbol\phi,\nabla\cdot\tau)_T
+\langle Q_b\boldsymbol\phi,\tau\cdot\bn\rangle_{\pT}\\
&=& -(\boldsymbol\phi,\nabla\cdot\tau)_T
+\langle \boldsymbol\phi,\tau\cdot\bn\rangle_{\pT}\\
&=&(\nabla \boldsymbol\phi,\tau)_T=(\Q_h\nabla\boldsymbol\phi,\tau)_T,
\end{eqnarray*}
which implies the identity (\ref{key1}).

Using (\ref{wd}) and  integration by parts, we have that for
any $w\in P_{k+1}(T)$
\begin{eqnarray*}
(\nabla_w \cdot Q_h\boldsymbol\phi,w)_T &=& -(Q_0\boldsymbol\phi,\nabla w)_T
+\langle Q_b\boldsymbol\phi\cdot\bn,w\rangle_{\pT}\\
&=& -(\boldsymbol\phi,\nabla w)_T
+\langle \boldsymbol\phi\cdot\bn,w\rangle_{\pT}\\
&=&(\nabla\cdot \boldsymbol\phi,w)_T=(\mathcal{Q}_h\nabla\cdot\boldsymbol\phi,w)_T,
\end{eqnarray*}
which proves (\ref{key1}).
\end{proof}

For any function $\varphi\in H^1(T)$, the following trace
inequality holds true (see \cite{wymix} for details):
\begin{equation}\label{trace}
\|\varphi\|_{e}^2 \leq C \left( h_T^{-1} \|\varphi\|_T^2 + h_T
\|\nabla \varphi\|_{T}^2\right).
\end{equation}

 We introduce two semi-norms $\3bar \bv\3bar$ and  $\|\bv\|_{1,h}$
   for any $\bv\in V_h$ as follows:
\begin{eqnarray}
\3bar \bv\3bar^2 &=& \sum_{T\in\T_h}(\nabla_w\bv,\nabla_w\bv)_T, \label{norm1}\\
\|\bv\|_{1,h}^2&=&\sum_{T\in \T_h}\|\nabla \bv_0\|_T^2+\sum_{T\in\T_h}h_T^{-1}\|\bv_0-\bv_b\|_{\pT}^2.\label{norm2}
\end{eqnarray}

It is easy to see that $\|\bv\|_{1,h}$ defines a norm in $V_h^0$.
Next we will show that $\3bar  \cdot \3bar$ also defines a norm in $V_h^0$ by proving the equivalence of $\3bar\cdot\3bar$ and $\|\cdot\|_{1,h}$ in $V_h$.

The following norm equivalence has been proved in  \cite{sf-wg-part-III} for each component of $\bv$,
\begin{equation}\label{happy}
C_1\|\bv\|_{1,h}\le \3bar \bv\3bar\le C_2 \|\bv\|_{1,h} \quad\forall \bv\in V_h.
\end{equation}


Unlike the traditional finite elements 
\cite{Li-M, Qin-P1,Qin-4,Zhang-M,Zhang-PS3,Zhang-PS2,Zhang-P6,Zhang-P2,Zhang-Z},
  the inf-sup condition for the
   weak Galerkin finite element is easily satisfied due to the large velocity space 
   with independent element boundary degrees of freedom.

\begin{lemma}\label{Lemma:inf-sup}
There exists a positive constant $\beta$ independent of $h$ such that for all $\rho\in W_h$,
\begin{equation}\label{inf-sup}
\sup_{\bv\in V_h}\frac{(\nabla_w\cdot\bv,\rho)}{\3bar\bv\3bar}\ge \beta
\|\rho\|.
\end{equation}
\end{lemma}

\begin{proof}
For any given $\rho\in W_h\subset L_0^2(\Omega)$, it is known
\cite{gr} that there exists a function
$\tilde\bv\in [H_0^1(\Omega)]^d$ such that
\begin{equation}\label{c-inf-sup}
\frac{(\nabla\cdot\tilde\bv,\rho)}{\|\tilde\bv\|_1}\ge C\|\rho\|,
\end{equation}
where $C>0$ is a constant independent of $h$.
Let $\bv=Q_h\tilde{\bv}=\{Q_0\tilde{\bv},Q_b\tilde{\bv}\}\in V_h$.
It follows from (\ref{happy}), (\ref{trace}) and $\tilde{\bv}\in [H_0^1(\Omega)]^d$,
\begin{eqnarray*}
\3bar \bv\3bar^2&\le& C\|\bv\|_{1,h}^2=C(\sum_{T\in \T_h}\|\nabla \bv_0\|_T^2+\sum_{T\in\T_h}h_T^{-1}\|\bv_0-\bv_b\|_{\pT}^2)\\
&\le&C (\sum_{T\in \T_h}\|\nabla Q_0\tilde{\bv}\|_T^2+\sum_{T\in\T_h}h_T^{-1}\|Q_0\tilde{\bv}-Q_b\tilde{\bv}\|_{\pT}^2)\\
&\le&C (\sum_{T\in \T_h}\|\nabla Q_0\tilde{\bv}\|_T^2+\sum_{T\in\T_h}h_T^{-1}\|Q_0\tilde{\bv}-\tilde{\bv}\|_{\pT}^2)\\
&\le& C\|\tilde{\bv}\|_1^2,
\end{eqnarray*}
which implies
\begin{equation}\label{m9}
\3bar\bv\3bar\le C\|\tilde{\bv}\|_1.
\end{equation}
It follows from (\ref{wd})  that
\begin{eqnarray}
(\nabla_w\cdot\bv,\;\rho)_{\T_h}&=&-(\bv_0,\;\nabla\rho)_{\T_h}+\l\bv_b,\rho\bn\r_{\partial\T_h}\nonumber\\
&=&-(Q_0\tilde\bv,\;\nabla\rho)_{\T_h}+\l Q_b\tilde\bv,\rho\bn\r_{\partial\T_h}\nonumber\\
&=&-(\tilde\bv,\;\nabla\rho)_{\T_h}+\l \tilde\bv,\rho\bn\r_{\partial\T_h}\nonumber\\
&=&(\nabla\cdot\tilde\bv,\;\rho)_{\T_h}.\label{m10}
\end{eqnarray}
Using  (\ref{m10}), (\ref{m9}) and (\ref{c-inf-sup}), we have
\begin{eqnarray*}
\frac{(\nabla_w\cdot\bv,\rho)} {\3bar\bv\3bar}=\frac{(\nabla\cdot\tilde{\bv},\rho)} {\3bar\bv\3bar}  &\ge &
\frac{(\nabla\cdot\tilde\bv,\rho)}{C\|\tilde\bv\|_1}\ge
\beta\|\rho\|,
\end{eqnarray*}
for a positive constant $\beta$. This completes the proof of the
lemma.
\end{proof}

\begin{lemma}
The weak Galerkin method (\ref{wg1})-(\ref{wg2}) has a unique solution.
\end{lemma}

\smallskip

\begin{proof}
It suffices to show that zero is the only solution of
(\ref{wg1})-(\ref{wg2}) if $\bbf=0$. To this end, let $\bbf=0$ and
take $\bv=\bu_h$ in (\ref{wg1}) and $w=p_h$ in (\ref{wg2}). By adding the
two resulting equations, we obtain
\[
(\nabla_w\bu_h,\ \nabla_w\bu_h)=0,
\]
which implies that $\nabla_w \bu_h=0$ on each element $T$. By (\ref{happy}), we have $\|\bu_h\|_{1,h}=0$ which implies that $\bu_h=0$.

Since $\bu_h=0$ and $\bbf=0$, the equation (\ref{wg1}) becomes $(\nabla\cdot\bv,\ p_h)=0$ for any $\bv\in V_h$. Then the inf-sup condition (\ref{inf-sup}) implies $p_h=0$. We have proved the lemma.
\end{proof}

\section{Error Equations}
In this section, we derive the equations that the errors satisfy.
Let $\be_h=Q_h\bu-\bu_h$ and $\varepsilon_h=\mathcal{Q}_hp-p_h$.

\begin{lemma}
The following error equations hold true for any $(\bv,w)\in V_h^0\times W_h$,
\begin{eqnarray}
(\nabla_w\be_h,\; \nabla_w\bv)-(\varepsilon_h,\;\nabla_w\cdot\bv)&=&\ell_1(\bu,\bv)+\ell_2(p,\bv),\label{ee1}\\
(\nabla_w\cdot\be_h,\ w)&=&0,\label{ee2}
\end{eqnarray}
where
\begin{eqnarray}
\ell_1(\bu,\ \bv)&=&(\Q_h\nabla\bu-\Pi_h\nabla\bu,\nabla_w\bv),\label{l1}\\
\ell_2(p,\;\bv)&=&\l \mathcal{Q}_hp-p, (\bv_0-\bv_b)\cdot\bn\r_{\partial\T_h}.\label{l3}
\end{eqnarray}
\end{lemma}

\begin{proof}
First, we test (\ref{moment}) by
$\bv_0$ with $\bv=\{\bv_0,\bv_b\}\in V_h^0$ to obtain
\begin{equation}\label{mm0}
-(\Delta\bu,\;\bv_0)+(\nabla p,\ \bv_0)=(\bbf,\; \bv_0).
\end{equation}
It follows from (\ref{key}) and (\ref{key1})
\begin{equation}\label{mm1}
-(\nabla\cdot\nabla\bu,\;\bv_0)=(\Pi_h\nabla\bu,\; \nabla_w\bv)=(\nabla_wQ_h\bu,\; \nabla_w\bv)-\ell_1(\bu,\bv).
\end{equation}
Using integration by parts and the fact $\l p, \bv_b\cdot\bn\r_{\partial\T_h}=0$,  we have
\begin{eqnarray*}
(\nabla p,\ \bv_0)&=& -(p,\nabla\cdot\bv_0)_{\T_h}+\l p, \bv_0\cdot\bn\r_{\partial\T_h}\\
&=& -(\mathcal{Q}_hp,\nabla\cdot\bv_0)_{\T_h}+\l p, (\bv_0-\bv_b)\cdot\bn\r_{\partial\T_h}\\
&=&(\nabla\mathcal{Q}_hp,\bv_0)_{\T_h}-\l \mathcal{Q}_hp, \bv_0\cdot\bn\r_{\partial\T_h}+\l p, (\bv_0-\bv_b)\cdot\bn\r_{\partial\T_h}\\
&=&-(\mathcal{Q}_hp,\nabla_w\cdot\bv)-\l \mathcal{Q}_hp, (\bv_0-\bv_b)\cdot\bn\r_{\partial\T_h}+\l p, (\bv_0-\bv_b)\cdot\bn\r_{\partial\T_h}\\
&=&-(\mathcal{Q}_hp,\nabla_w\cdot\bv)-\ell_2(p,\bv),
\end{eqnarray*}
which implies
\begin{equation}\label{mm2}
(\nabla p,\ \bv_0)=-(\mathcal{Q}_hp,\nabla_w\cdot\bv)-\ell_2(p,\bv).
\end{equation}
Substituting (\ref{mm1}) and (\ref{mm2}) into (\ref{mm0}) gives
\begin{equation}\label{mm3}
(\nabla_wQ_h\bu,\nabla_w\bv)-(\mathcal{Q}_hp,\nabla_w\cdot\bv)=(\bbf,\bv_0)+\ell_1(\bu,\bv)+\ell_2(p,\bv).
\end{equation}
The difference of (\ref{mm3}) and (\ref{wg1}) implies
\begin{equation}\label{mm4}
(\nabla_w\be_h,\nabla_w\bv)-(\varepsilon_h,\nabla_w\cdot\bv)=\ell_1(\bu,\bv)+\ell_2(p,\bv)\quad\forall\bv\in V_h^0.
\end{equation}
Testing equation (\ref{cont}) by $w\in W_h$ and using (\ref{key2}) give
\begin{equation}\label{mm5}
(\nabla\cdot\bu,\ w)=(\mathcal{Q}_h\nabla\cdot\bu,\ w)=(\nabla_w\cdot Q_h\bu,\ w)=0.
\end{equation}
The difference of (\ref{mm5}) and (\ref{wg2}) implies (\ref{ee2}). We have proved the lemma.
\end{proof}

\section{Error Estimates in Energy Norm}\label{Section:error-analysis}
In this section, we establish order two superconvergence
for the velocity approximation $\bu_h$ in $\3bar\cdot\3bar$ norm and optimal order error estimate for the pressure approximation $p_h$ in
the standard $L^2$ norm.

\begin{lemma}
Let $\bu\in [H^{k+3}(\Omega)]^d$ and $p\in H^{k+2}(\Omega)$ and
$\bv\in V_h$. Then, the following estimates hold true
\begin{eqnarray}
|\ell_1(\bu,\ \bv)|&\le& Ch^{k+2}|\bu|_{k+3}\3bar \bv\3bar,\label{mmm2}\\
|\ell_2(p,\ \bv)|&\le& Ch^{k+2}|p|_{k+2}\3bar \bv\3bar.\label{mmm3}
\end{eqnarray}
\end{lemma}

\begin{proof}
Using the Cauchy-Schwarz inequality and the definitions of $\Q_h$ and $\Pi_h$, we have
\begin{eqnarray*}
|\ell_1(\bu,\ \bv)|&=&|(\Q_h\nabla\bu-\Pi_h\nabla\bu,\nabla_w\bv)|\\
&=&|(\Q_h\nabla\bu-\nabla\bu+\nabla\bu-\Pi_h\nabla\bu,\nabla_w\bv)|\\
&\le & Ch^{k+2}|\bu|_{k+3}\3bar \bv\3bar.
\end{eqnarray*}
It follows from (\ref{trace}) and (\ref{happy})
\begin{eqnarray*}
|\ell_2(p,\;\bv)|&=&|\l \mathcal{Q}_hp-p, (\bv_0-\bv_b)\cdot\bn\r_{\partial\T_h}|\\
&\le & C \sum_{T\in\T_h}\|\mathcal{Q}_hp-p\|_{\pT}
\|\bv_0-\bv_b\|_\pT\nonumber\\
&\le & C \left(\sum_{T\in\T_h}h_T\|\mathcal{Q}_hp-p\|_{\pT}^2\right)^{\frac12}
\left(\sum_{e\in\E_h}h_T^{-1}\|\bv_0-\bv_b\|_e^2\right)^{\frac12}\\
&\le& Ch^{k+2}|p|_{k+2}\3bar \bv\3bar.
\end{eqnarray*}
We have proved the lemma.
\end{proof}

\begin{theorem}\label{h1-bd}
Let  $(\bu_h,p_h)\in
V_h^0\times W_h$ be the solution of (\ref{wg1})-(\ref{wg2}). Then, we have
\begin{eqnarray}
\3bar Q_h\bu-\bu_h\3bar &\le& Ch^{k+2}(|\bu|_{k+3}+|p|_{k+2}),\label{errv}\\
\|\mathcal{Q}_hp-p_h\|&\le& Ch^{k+2}(|\bu|_{k+3}+|p|_{k+2}).\label{errp}
\end{eqnarray}
\end{theorem}

\smallskip

\begin{proof}
By letting $\bv=\be_h$ in (\ref{ee1}) and $w=\varepsilon_h$ in
(\ref{ee2}) and using the equation (\ref{ee2}), we have
\begin{eqnarray}
\3bar \be_h\3bar^2&=&|\ell_1(\bu,\be_h)+\ell_2(p,\be_h)|.\label{main}
\end{eqnarray}
It then follows from (\ref{mmm2}) and (\ref{mmm3}) that
\begin{equation}\label{b-u}
\3bar \be_h\3bar^2 \le Ch^{k+2}(|\bu|_{k+3}+|p|_{k+2})\3bar \be_h\3bar,
\end{equation}
which implies (\ref{errv}). To estimate
$\|\varepsilon_h\|$, we have from (\ref{ee1}) that
\[
(\varepsilon_h, \nabla\cdot\bv)=(\nabla_w\be_h,\nabla_w\bv)-\ell_1(\bu, \bv)-\ell_2(p, \bv).
\]
Using (\ref{b-u}), (\ref{mmm2}) and
(\ref{mmm3}), we arrive at
\[
|(\varepsilon_h, \nabla\cdot\bv)|\le Ch^{k+2}(|\bu|_{k+3}+|p|_{k+2})\3bar\bv\3bar.
\]
Combining the above estimate with the {\em inf-sup} condition
(\ref{inf-sup}) gives
\[
\|\varepsilon_h\|\le Ch^{k+2}(|\bu|_{k+3}+|p|_{k+2}),
\]
which yields the desired estimate (\ref{errp}).
\end{proof}

\section{Error Estimates in $L^2$ Norm}

In this section, order two superconvergence for velocity in the $L^2$ norm is obtained by duality argument. Recall that $\be_h=\{\be_0,\be_b\}=Q_h\bu-\bu_h$ and $\epsilon_h=\mathcal{Q}_hp-p_h$. Consider the dual problem: seeking $(\bpsi,\xi)$ satisfying 
\begin{eqnarray}
-\Delta\bpsi+\nabla \xi&=\be_0 &\quad \mbox{in}\;\Omega,\label{dual-m}\\
\nabla\cdot\bpsi&=0 &\quad\mbox{in}\;\Omega,\label{dual-c}\\
\bpsi&= 0 &\quad\mbox{on}\;\partial\Omega.\label{dual-bc}
\end{eqnarray}
Assume that the dual problem (\ref{dual-m})-(\ref{dual-bc}) satisfy the following regularity assumption:
\begin{eqnarray}
\|\bpsi\|_{2}+\|\xi\|_1&\le& C\|\be_0\|.\label{reg}
\end{eqnarray}
We need the following lemma first.

\begin{lemma}
For any $\bv\in V_h^0$ and $w\in W_h$, the following equations hold true,
\begin{eqnarray}
(\nabla_w Q_h\bpsi,\; \nabla_w\bv)-(\mathcal{Q}_h\xi,\;\nabla_w\cdot\bv)&=&(\be_0,\bv_0)+\ell_3(\bpsi,\bv)+\ell_2(\xi,\bv),\label{wd1}\\
(\nabla_w\cdot Q_h\bpsi,\ w)&=&0,\label{wd2}
\end{eqnarray}
where
\begin{eqnarray*}
\ell_3(\bpsi,\ \bv)&=&\langle (\nabla\bpsi-\Q_h\nabla\bpsi)\cdot\bn,\;\bv_0-\bv_b\rangle_{\partial\T_h},\\
\ell_2(\xi,\;\bv)&=&\l \mathcal{Q}_h\xi-\xi, (\bv_0-\bv_b)\cdot\bn\r_{\partial\T_h}.
\end{eqnarray*}
\end{lemma}

\begin{proof}
Testing (\ref{dual-m}) by
$\bv_0$ with $\bv=\{\bv_0,\bv_b\}\in V_h^0$ gives
\begin{equation}\label{d1}
-(\Delta\bpsi,\;\bv_0)+(\nabla \xi,\ \bv_0)=(\be_0,\; \bv_0).
\end{equation}
It follows from integration by parts and the fact $\langle\nabla \bpsi\cdot\bn,\bv_b\rangle_{\partial\T_h}=0$
\begin{equation}\label{d2}
-(\Delta\bpsi,\;\bv_0)=(\nabla\bpsi,\nabla \bv_0)_{\T_h}- \langle
\nabla \bpsi\cdot\bn,\bv_0-\bv_b\rangle_{\partial\T_h}.
\end{equation}
By integration by parts, (\ref{wg}) and (\ref{key1})
\begin{eqnarray}
(\nabla \bpsi,\nabla \bv_0)_{\T_h}&=&(\Q_h\nabla\bpsi,\nabla \bv_0)_{\T_h}\nonumber\\
&=&-(\bv_0,\nabla\cdot (\Q_h\nabla \bpsi))_{\T_h}+\langle \bv_0, \Q_h\nabla \bpsi\cdot\bn\rangle_{\partial\T_h}\nonumber\\
&=&(\Q_h\nabla \bpsi, \nabla_w \bv)+\langle \bv_0-\bv_b,\Q_h\nabla \psi\cdot\bn\rangle_{\partial\T_h}\nonumber\\
&=&( \nabla_w Q_h\bpsi, \nabla_w \bv)+\langle \bv_0-\bv_b,\Q_h\nabla \bpsi\cdot\bn\rangle_{\partial\T_h}.\label{d3}
\end{eqnarray}
Combining (\ref{d2}) and (\ref{d3}) gives
\begin{eqnarray}
-(\Delta\bpsi,\;\bv_0)&=&( \nabla_w Q_h\bpsi, \nabla_w \bv)-\ell_3(\bpsi,\bv).\label{d4}
\end{eqnarray}
Similar to the derivation of (\ref{mm2}), we obtain
\begin{equation}\label{d5}
(\nabla \xi,\ \bv_0)=-(\mathcal{Q}_h\xi,\nabla_w\cdot\bv)-\ell_2(\xi,\bv).
\end{equation}
Combining (\ref{d4}) and (\ref{d5}) with (\ref{d1}) yields (\ref{wd1}).
Testing equation (\ref{dual-c}) by $w\in W_h$ and using (\ref{key2}) give
\begin{equation}\label{d6}
(\nabla\cdot\bpsi,\ w)=(\mathcal{Q}_h\nabla\cdot\bpsi,\ w)=(\nabla_w\cdot Q_h\bpsi,\ w)=0,
\end{equation}
which implies (\ref{wd2}) and we have proved the lemma.
\end{proof}

\smallskip

By the same argument as (\ref{d4}), (\ref{ee1}) has another form as
\begin{eqnarray}
(\nabla_w \be_h,\; \nabla_w\bv)-(\epsilon_h,\;\nabla_w\cdot\bv)&=&\ell_3(\bu,\bv)+\ell_2(p,\bv).\label{eee1}
\end{eqnarray}

\begin{theorem}
Let  $(\bu_h,p_h)\in V_h^0\times W_h$ be the solution of
(\ref{wg1})-(\ref{wg2}). Assume that (\ref{reg}) holds true.  Then, we have
\begin{equation}\label{l2err}
\|Q_0\bu-\bu_0\|\le Ch^{k+3}(|\bu|_{k+3}+|p|_{k+2}).
\end{equation}
\end{theorem}

\begin{proof}
Letting $\bv=\be_h$ in (\ref{wd1}) yields
\begin{eqnarray}
\|\be_0\|^2=(\nabla_w Q_h\bpsi,\; \nabla_w\be_h)-(Q_h\xi,\;\nabla_w\cdot\be_h)-\ell_3(\bpsi,\be_h)-\ell_2(\xi,\be_h).\label{d7}
\end{eqnarray}
Using the fact $(Q_h\xi,\;\nabla_w\cdot\be_h)=0$, (\ref{d7}) becomes
\begin{eqnarray}
\|\be_h\|^2=(\nabla_w Q_h\bpsi,\; \nabla_w\be_h)-\ell_3(\bpsi,\be_h)-\ell_2(\xi,\be_h).\label{d8}
\end{eqnarray}
With $\bv=Q_h\bpsi$, (\ref{eee1}) becomes
\begin{eqnarray}
(\nabla_w \be_h,\; \nabla_w Q_h\bpsi)-(\epsilon_h,\;\nabla_w\cdot Q_h\bpsi)&=&\ell_3(\bu,Q_h\bpsi)+\ell_2(p,Q_h\bpsi).\label{d9}
\end{eqnarray}
Using (\ref{wd2}), we have $(\epsilon_h,\;\nabla_w\cdot Q_h\bpsi)=0$. Then (\ref{d9}) becomes
\begin{eqnarray}
(\nabla_w \be_h,\; \nabla_w Q_h\bpsi)&=&\ell_3(\bu,Q_h\bpsi)+\ell_2(p,Q_h\bpsi).\label{d10}
\end{eqnarray}
Combining (\ref{d8}) and (\ref{d10}), we have
\begin{eqnarray}
\|\be_h\|^2&=&\ell_3(\bu,Q_h\bpsi)+\ell_2(p,Q_h\bpsi)-\ell_3(\bpsi,\be_h)-\ell_2(\xi,\be_h).\label{d11}
\end{eqnarray}
 Using the Cauchy-Schwarz inequality, the trace inequality (\ref{trace}) and the definition of $\Q_h$,
we arrive at
\begin{eqnarray}
|\ell_3(\bu,Q_h\bpsi)|&\le&\left| \langle (\nabla \bu-\Q_h\nabla
\bu)\cdot\bn,\;
Q_0\bpsi-Q_b\bpsi\rangle_{\pT_h} \right|\nonumber\\
&\le& \left(\sum_{T\in\T_h}\|\nabla \bu-\Q_h\nabla \bu\|^2_\pT\right)^{1/2}
\left(\sum_{T\in\T_h}\|Q_0\bpsi-\bpsi\|^2_\pT\right)^{1/2}\nonumber \\
&\le&  Ch^{k+3}|\bu|_{k+3}|\bpsi|_2.\label{e1}
\end{eqnarray}
Similarly, we have
\begin{eqnarray}
|\ell_2(p, Q_h\bpsi)|&\le&\left| \langle \mathcal{Q}_hp-p,\;
(Q_0\bpsi-Q_b\bpsi)\cdot\bn\rangle_{\pT_h} \right|\nonumber\\
&\le& C \left(\sum_{T\in\T_h}\|\mathcal{Q}_hp-p\;\|^2_\pT\right)^{1/2}
\left(\sum_{T\in\T_h}\|Q_0\bpsi-\bpsi\|^2_\pT\right)^{1/2}\nonumber \\
&\le&  Ch^{k+3}|p|_{k+2}|\bpsi|_2.\label{e2}
\end{eqnarray}
It follows from the Cauchy-Schwarz inequality, the trace inequality (\ref{happy}) and (\ref{errv}), 
\begin{eqnarray}
|\ell_3(\bpsi,\be_h)|&\le&\left| \langle (\nabla \bpsi-\Q_h\nabla
\bpsi)\cdot\bn,\;
\be_0-\be_b\rangle_{\pT_h} \right|\nonumber\\
&\le& \left(\sum_{T\in\T_h}h_T\|\nabla \bpsi-\Q_h\nabla \bpsi\|^2_\pT\right)^{1/2}
\left(\sum_{T\in\T_h}h_T^{-1}\|\be_0-\be_b\|^2_\pT\right)^{1/2}\nonumber \\
&\le&  Ch|\bpsi|_2\3bar\be_h\3bar\nonumber\\
&\le&Ch^{k+3}(|\bu|_{k+3}+|p|_{k+2})|\bpsi|_2.\label{e3}
\end{eqnarray}
Similarly, 
\begin{eqnarray}
|\ell_2(\xi, \be_h)|&\le&\left| \langle \mathcal{Q}_h\xi-\xi,\;
(\be_0-\be_b)\cdot\bn\rangle_{\pT_h} \right|\nonumber\\
&\le& \left(\sum_{T\in\T_h}h_T\|\mathcal{Q}_h\xi-\xi\;\|^2_\pT\right)^{1/2}
\left(\sum_{T\in\T_h}h_T^{-1}\|\be_0-\be_b\|^2_\pT\right)^{1/2}\nonumber \\
&\le&  Ch^{k+3}(|\bu|_{k+3}+|p|_{k+2})|\xi|_1.\label{e4}
\end{eqnarray}
Combining all the estimates  above with (\ref{d11}) yields
$$
\|\be_h\|^2 \leq C h^{k+3}(|\bu|_{k+3}+|p|_{k+2})(\|\bpsi\|_2+\|\xi\|_1).
$$
The estimate (\ref{l2err}) follows from the above inequality and
the regularity assumption (\ref{reg}). We have completed the proof.
\end{proof}

\section{Numerical Experiments}\label{Section:numerical-experiments}

\begin{figure}[htb]\begin{center}\setlength\unitlength{1.5in}
    \begin{picture}(3.2,1.4)
 \put(0,0){\includegraphics[width=1.5in]{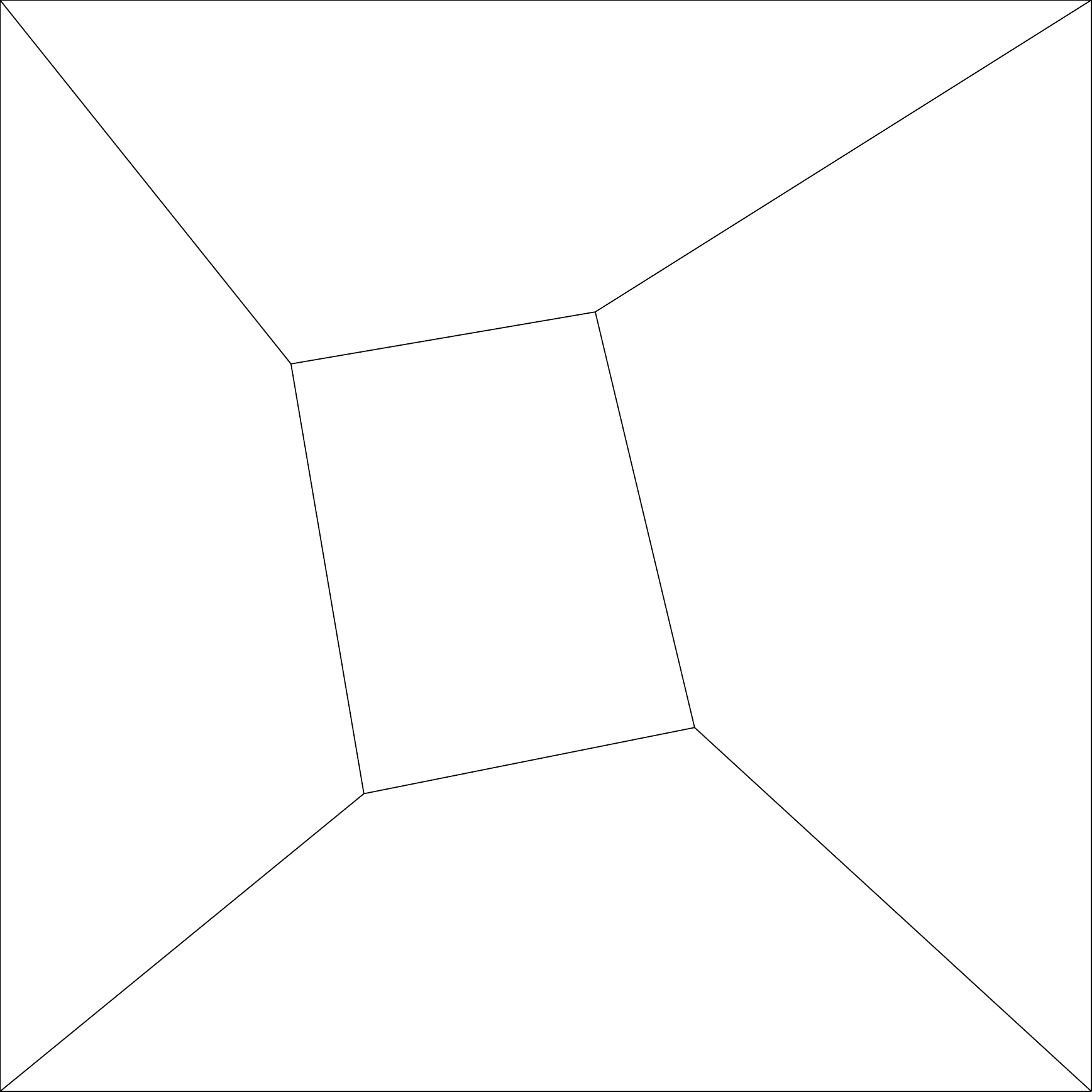}}
 \put(1.1,0){\includegraphics[width=1.5in]{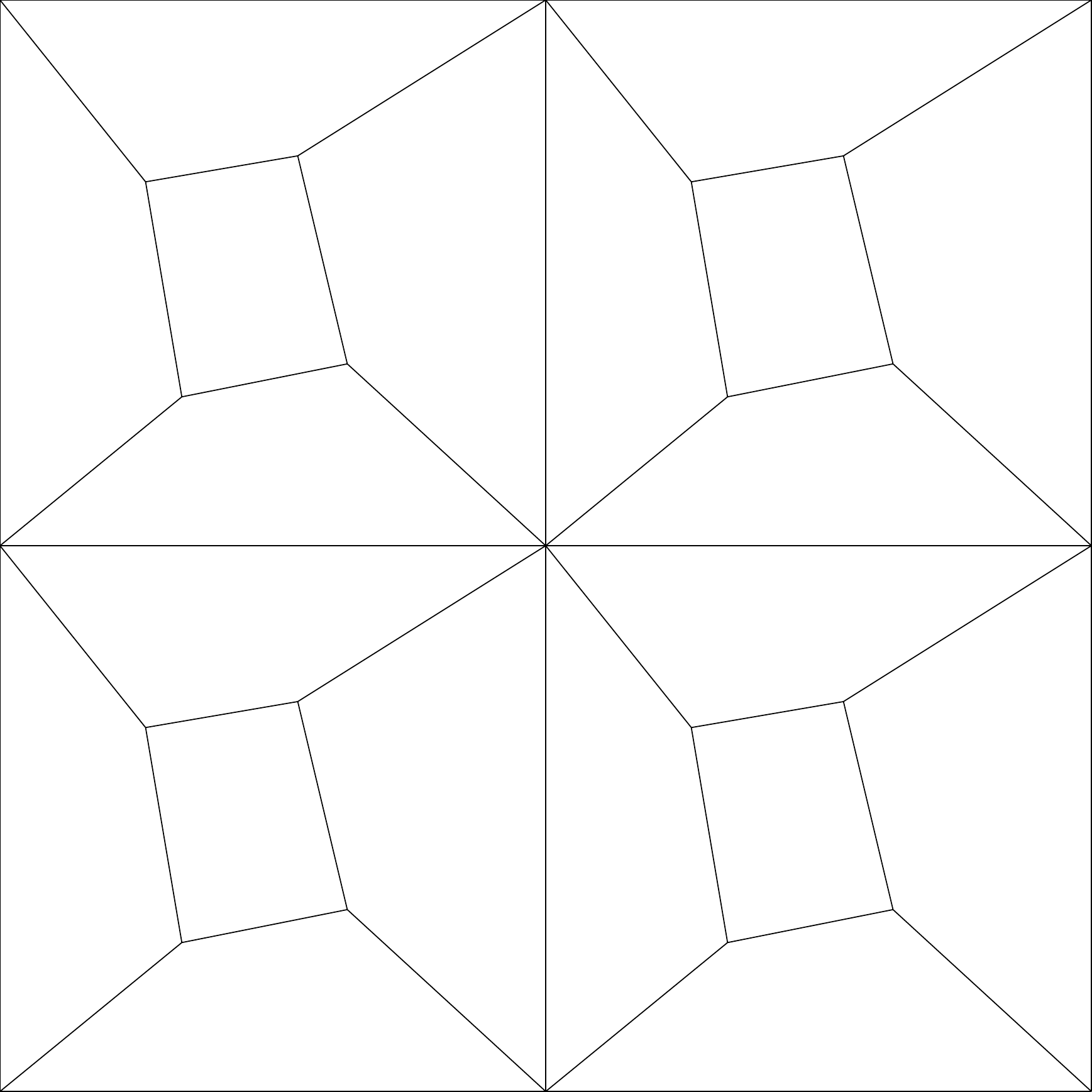}}
 \put(2.2,0){\includegraphics[width=1.5in]{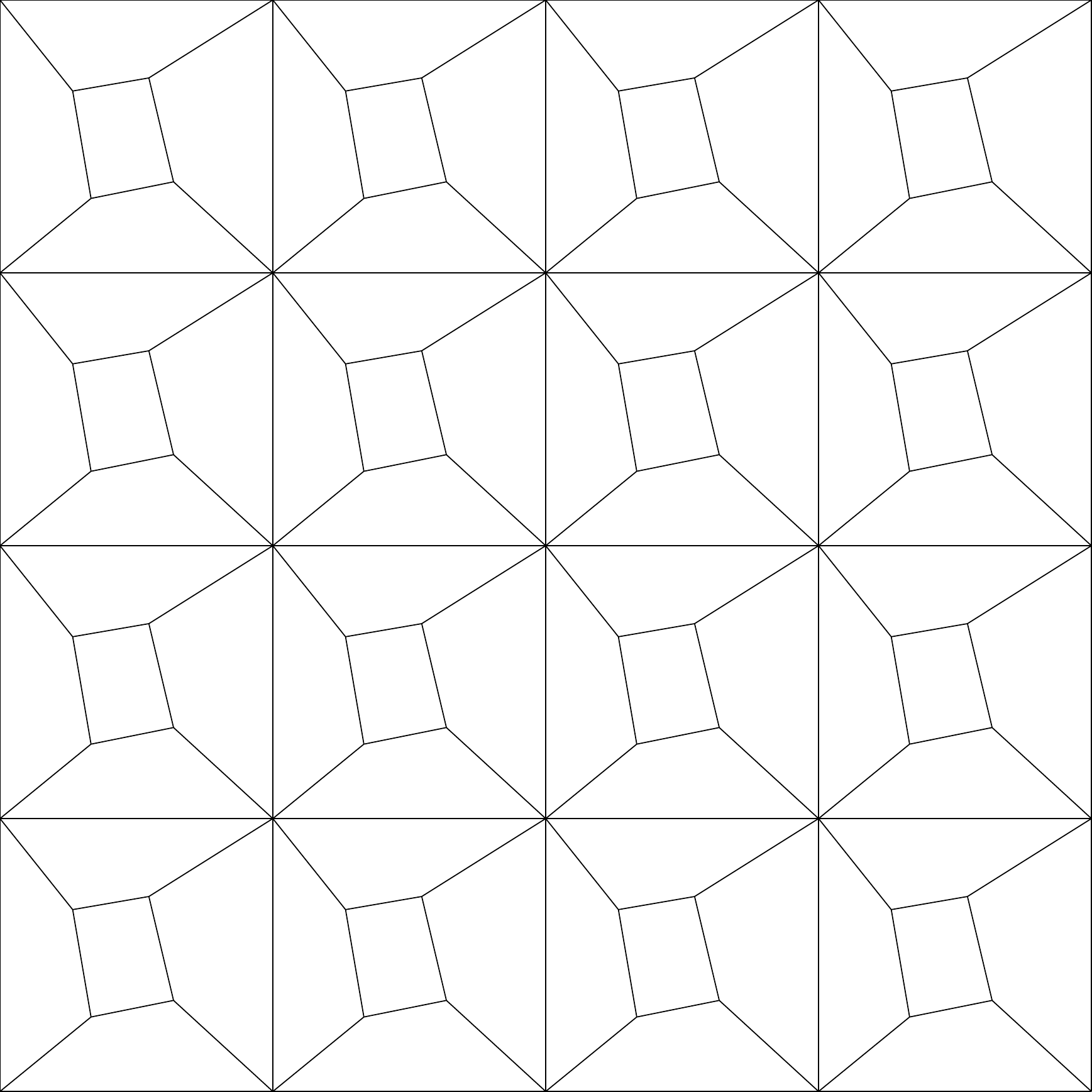}}
    \end{picture}
\caption{ The first three quadrilateral grids
   for the computation of Table \ref{t1}.  } \label{grid-t}
\end{center}
\end{figure}

\begin{table}[ht]
  \centering   \renewcommand{\arraystretch}{1.05}
  \caption{Error profiles and convergence rates for solution \eqref{s1}
     on quadrilateral grids shown in Figure \ref{grid-t}. }
\label{t1}
\begin{tabular}{c|cc|cc|cc}
\hline
Grid & $\|Q_h \bu- \bu_h \|_0 $  &rate &  $\3bar Q_h \bu- \bu_h \3bar $ &rate
  &  $\|p - p_h \|_0 $ &rate   \\
\hline
 &\multicolumn{6}{c}{by the $P_0^2$-$P_1^2$-$P_1$ WG finite element } \\ \hline
 4&   0.2179E-01 &  1.81&   0.2970E+00 &  1.92&   0.2118E+00 &  1.95 \\
 5&   0.5640E-02 &  1.95&   0.7565E-01 &  1.97&   0.5350E-01 &  1.99 \\
 6&   0.1422E-02 &  1.99&   0.1905E-01 &  1.99&   0.1347E-01 &  1.99 \\
\hline
 &\multicolumn{6}{c}{by the $P_1^2$-$P_2^2$-$P_2$ WG finite element } \\ \hline
 4&   0.3051E-03 &  3.95&   0.3440E-01 &  3.02&   0.9223E-02 &  2.95 \\
 5&   0.1964E-04 &  3.96&   0.4313E-02 &  3.00&   0.1209E-02 &  2.93 \\
 6&   0.1248E-05 &  3.98&   0.5421E-03 &  2.99&   0.1555E-03 &  2.96 \\
\hline
 &\multicolumn{6}{c}{by the $P_2^2$-$P_3^2$-$P_3$ WG finite element } \\ \hline
 3&   0.8289E-03 &  5.12&   0.8054E-01 &  4.12&   0.5896E-02 &  4.32 \\
 4&   0.2507E-04 &  5.05&   0.4871E-02 &  4.05&   0.3609E-03 &  4.03 \\
 5&   0.7763E-06 &  5.01&   0.3018E-03 &  4.01&   0.2277E-04 &  3.99 \\
 \hline
 &\multicolumn{6}{c}{by the $P_3^2$-$P_3^2$-$P_3$ WG finite element } \\ \hline
 2&   0.6018E-02 &  6.29&   0.3910E+00 &  5.28&   0.1249E-01 &  5.72 \\
 3&   0.8806E-04 &  6.09&   0.1146E-01 &  5.09&   0.2933E-03 &  5.41 \\
 4&   0.1352E-05 &  6.03&   0.3526E-03 &  5.02&   0.8304E-05 &  5.14 \\
 \hline
\end{tabular}%
\end{table}%

Consider problem  \eqref{moment}--\eqref{bc} with $\Omega=(0,1)^2$.
The source term $\bf$ and the boundary value $\bg$ are chosen so that the exact solution is
\an{\label{s1}
    \bu&=\p{g_y \\ -g_x}, \quad
      p =g_{xy}, \quad\hbox{where } \ g=2^4(x-x^2)^2(y-y^2)^2.
}
In this example, we use quadrilateral grids shown in Figure \ref{grid-t}.
In Table \ref{t1}, we list the errors and the orders of convergence.
We can see that two-order superconvergence is achieved for the
   velocity in $L^2$-norm and $H^1$-like norm.
The pressure converges at the optimal order.

\begin{figure}[htb]\begin{center}\setlength\unitlength{1.5in}
    \begin{picture}(3.2,1.4)
 \put(0,0){\includegraphics[width=1.5in]{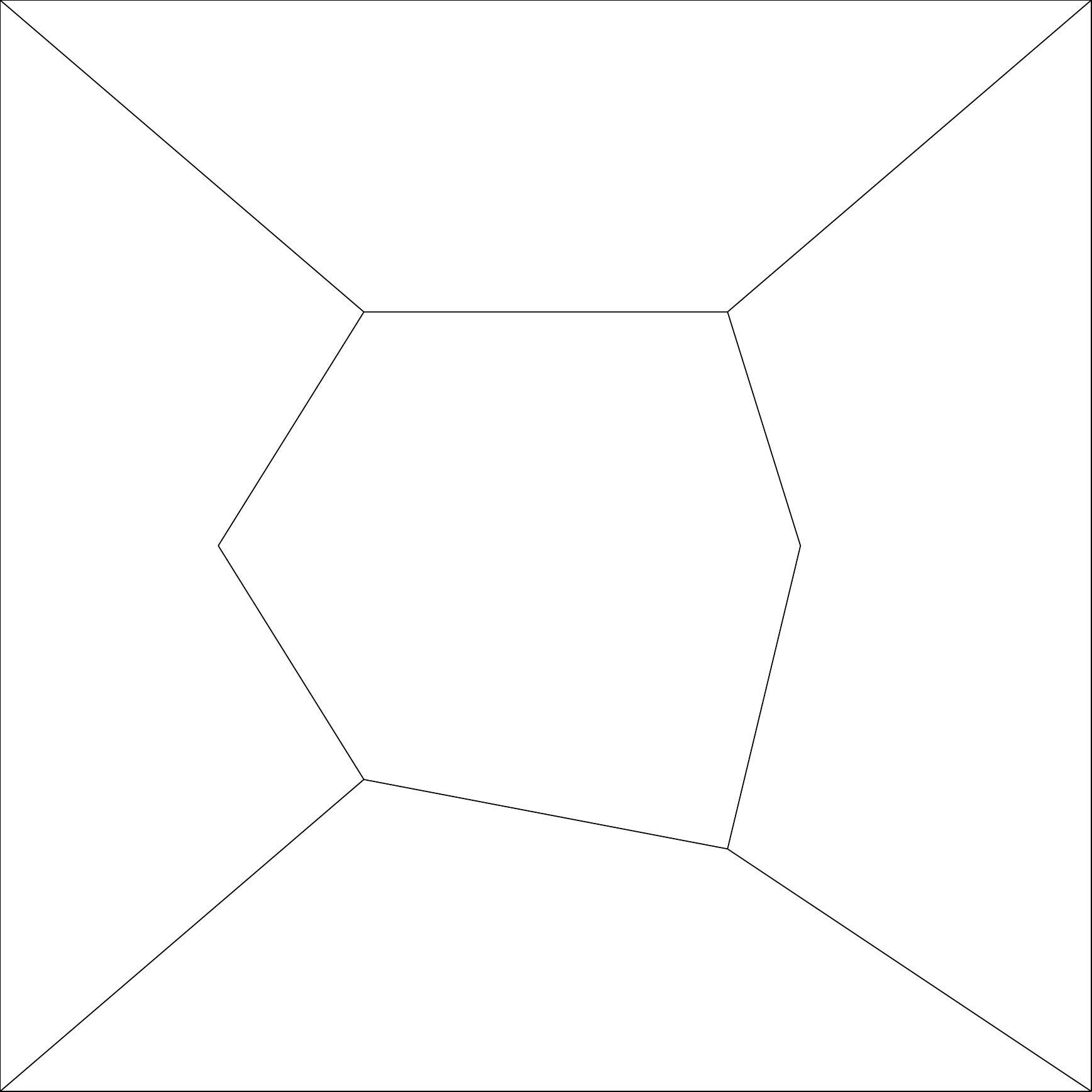}}
 \put(1.1,0){\includegraphics[width=1.5in]{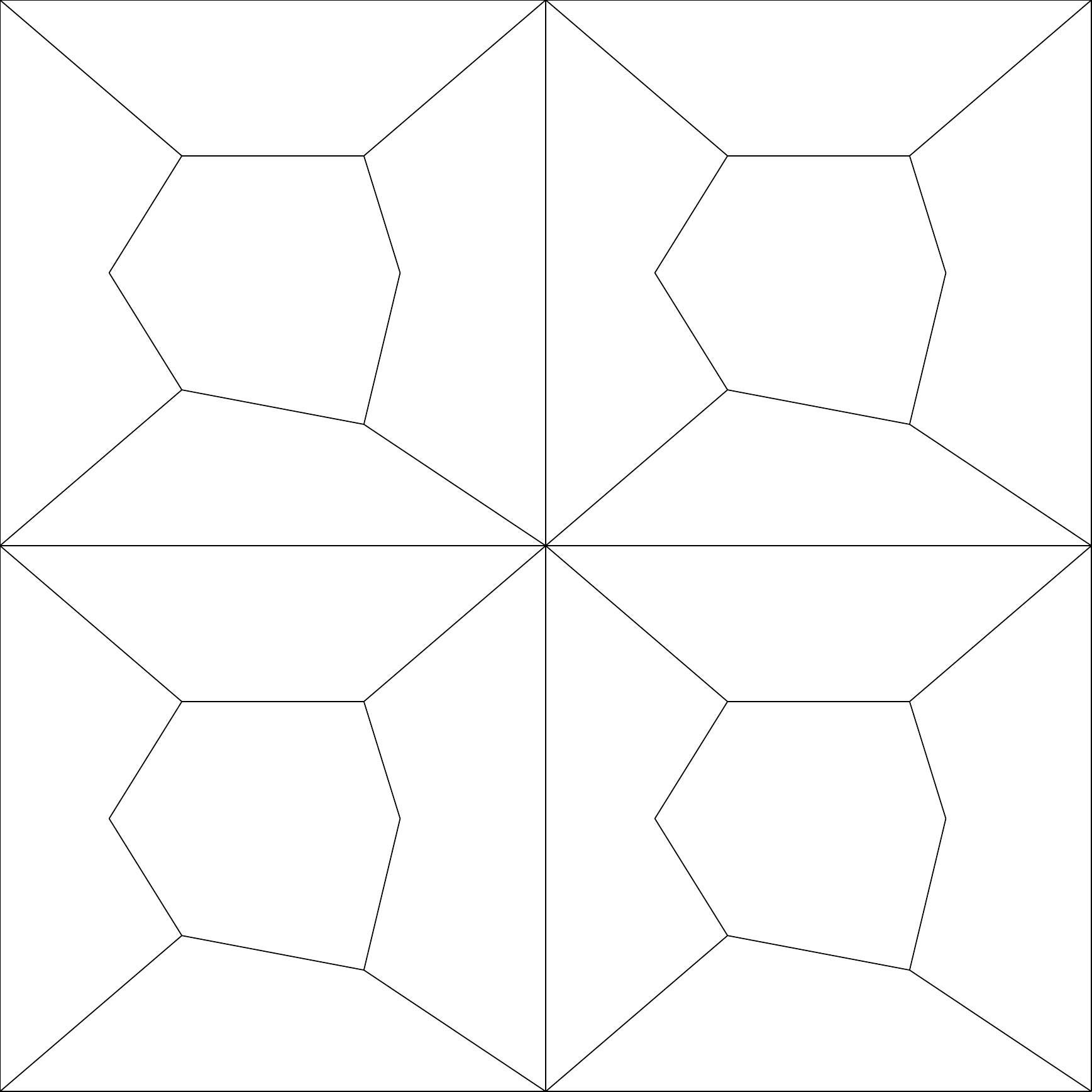}}
 \put(2.2,0){\includegraphics[width=1.5in]{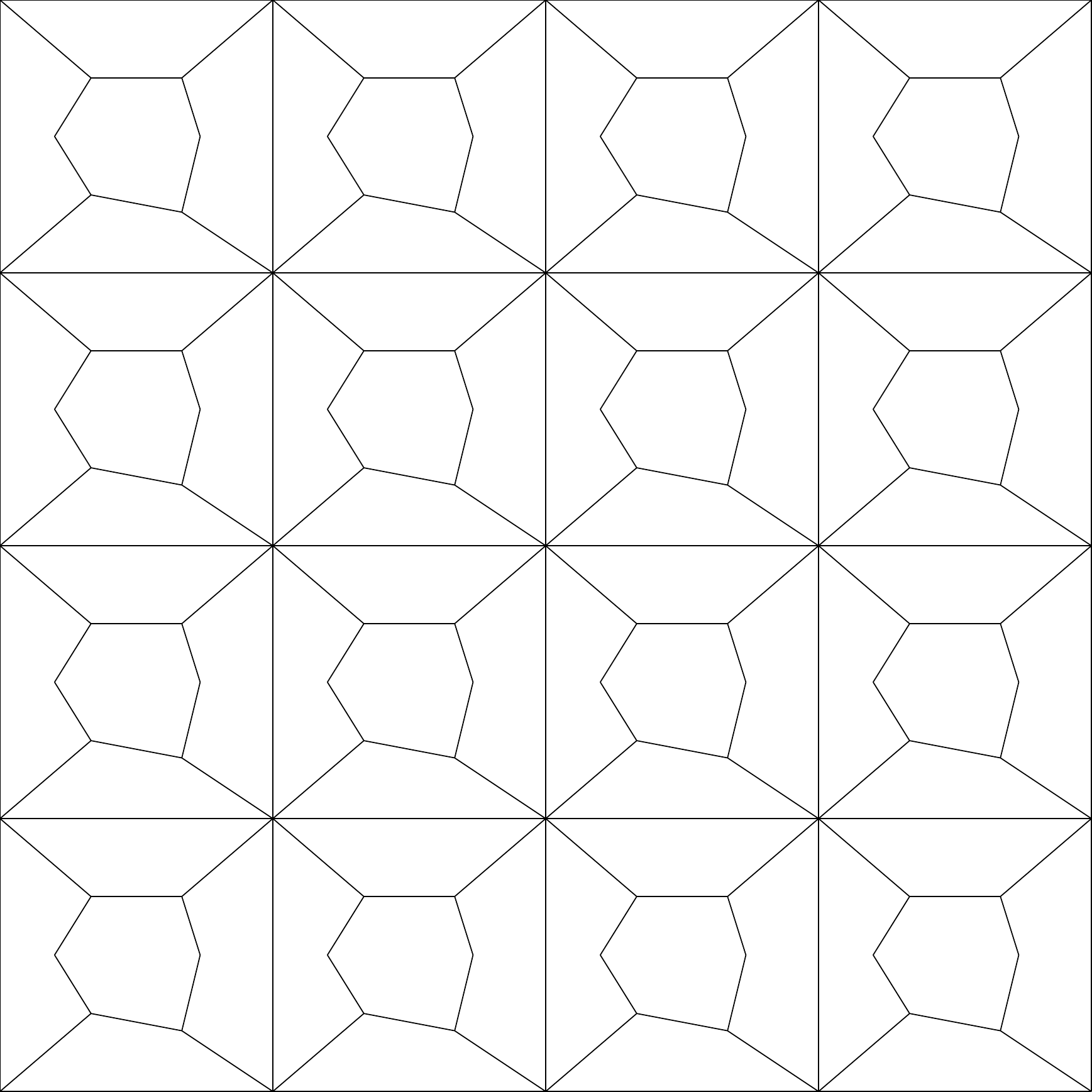}}
    \end{picture}
\caption{ The first three polygonal grids
   for the computation of Table \ref{t2}.  } \label{grid-6}
\end{center}
\end{figure}

We solve above problem \eqref{s1} again,  on polygonal grids, consisting of
     quadrilaterals, pentagons and hexagons,
      shown in Figure \ref{grid-6}.
In Table \ref{t2}, we list the errors and the orders of convergence.
The computational results match the theoretic order of convergence, in all cases.

\begin{table}[ht]
  \centering   \renewcommand{\arraystretch}{1.05}
  \caption{Error profiles for solution \eqref{s1} on polygonal grids shown
       in Figure \ref{grid-6}. }
\label{t2}
\begin{tabular}{c|cc|cc|cc}
\hline
Grid & $\|Q_h \bu- \bu_h \|_0 $  &rate &  $\3bar Q_h \bu- \bu_h \3bar $ &rate
  &  $\|p - p_h \|_0 $ &rate   \\
\hline
 &\multicolumn{6}{c}{by the $P_0^2$-$P_1^2$-$P_1$ WG finite element } \\ \hline
 4&   0.2202E-01 &  1.80&   0.2885E+00 &  1.91&   0.2138E+00 &  1.97 \\
 5&   0.5715E-02 &  1.95&   0.7374E-01 &  1.97&   0.5376E-01 &  1.99 \\
 6&   0.1442E-02 &  1.99&   0.1861E-01 &  1.99&   0.1351E-01 &  1.99 \\
\hline
 &\multicolumn{6}{c}{by the $P_1^2$-$P_2^2$-$P_2$ WG finite element } \\ \hline
 4&   0.2512E-03 &  3.86&   0.2922E-01 &  2.94&   0.8673E-02 &  2.93 \\
 5&   0.1661E-04 &  3.92&   0.3737E-02 &  2.97&   0.1147E-02 &  2.92 \\
 6&   0.1066E-05 &  3.96&   0.4735E-03 &  2.98&   0.1481E-03 &  2.95 \\
\hline
 &\multicolumn{6}{c}{by the $P_2^2$-$P_3^2$-$P_3$ WG finite element } \\ \hline
 3&   0.5373E-03 &  5.10&   0.5945E-01 &  4.16&   0.5018E-02 &  4.11 \\
 4&   0.1639E-04 &  5.04&   0.3567E-02 &  4.06&   0.3219E-03 &  3.96 \\
 5&   0.5101E-06 &  5.01&   0.2208E-03 &  4.01&   0.2063E-04 &  3.96 \\

 \hline
 &\multicolumn{6}{c}{by the $P_3^2$-$P_3^2$-$P_3$ WG finite element } \\ \hline
 2&   0.3384E-02 &  6.26&   0.2770E+00 &  5.30&   0.8473E-02 &  5.54 \\
 3&   0.4985E-04 &  6.08&   0.8069E-02 &  5.10&   0.2273E-03 &  5.22 \\
 4&   0.7855E-06 &  5.99&   0.2525E-03 &  5.00&   0.6876E-05 &  5.05 \\ \hline
\end{tabular}%
\end{table}%

\begin{figure}[ht]
\begin{center}
 \setlength\unitlength{1pt}
    \begin{picture}(320,118)(0,3)
    \put(0,0){\begin{picture}(110,110)(0,0)
       \multiput(0,0)(80,0){2}{\line(0,1){80}}  \multiput(0,0)(0,80){2}{\line(1,0){80}}
       \multiput(0,80)(80,0){2}{\line(1,1){20}} \multiput(0,80)(20,20){2}{\line(1,0){80}}
       \multiput(80,0)(0,80){2}{\line(1,1){20}}  \multiput(80,0)(20,20){2}{\line(0,1){80}}
    \put(80,0){\line(-1,1){80}}
      \end{picture}}
    \put(110,0){\begin{picture}(110,110)(0,0)
       \multiput(0,0)(40,0){3}{\line(0,1){80}}  \multiput(0,0)(0,40){3}{\line(1,0){80}}
       \multiput(0,80)(40,0){3}{\line(1,1){20}} \multiput(0,80)(10,10){3}{\line(1,0){80}}
       \multiput(80,0)(0,40){3}{\line(1,1){20}}  \multiput(80,0)(10,10){3}{\line(0,1){80}}
    \put(80,0){\line(-1,1){80}}
       \multiput(40,0)(40,40){2}{\line(-1,1){40}}
      \end{picture}}
    \put(220,0){\begin{picture}(110,110)(0,0)
       \multiput(0,0)(20,0){5}{\line(0,1){80}}  \multiput(0,0)(0,20){5}{\line(1,0){80}}
       \multiput(0,80)(20,0){5}{\line(1,1){20}} \multiput(0,80)(5,5){5}{\line(1,0){80}}
       \multiput(80,0)(0,20){5}{\line(1,1){20}}  \multiput(80,0)(5,5){5}{\line(0,1){80}}
    \put(80,0){\line(-1,1){80}}
       \multiput(40,0)(40,40){2}{\line(-1,1){40}}

       \multiput(20,0)(60,60){2}{\line(-1,1){20}}   \multiput(60,0)(20,20){2}{\line(-1,1){60}}
      \end{picture}}

    \end{picture}
    \end{center}
\caption{  The first three levels of wedge grids used in Table \ref{t3}. }
\label{grid3}
\end{figure}
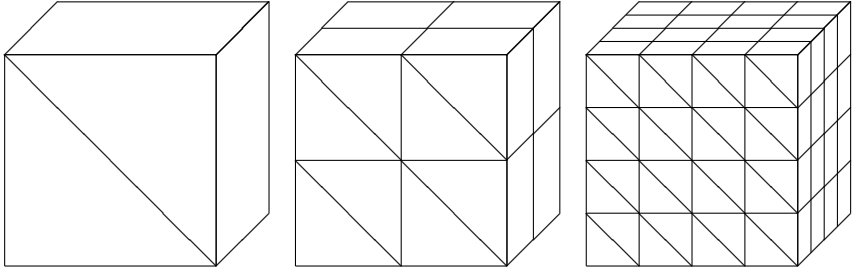

Finally we compute a 3D problem \eqref{moment}--\eqref{bc} with $\Omega=(0,1)^3$.
The source term $\bf$ and the boundary value $\bg$ are chosen so that the exact solution is
\an{\label{s2}
    \bu &=\p{-g_y\\g_x+g_z\\-g_y  }, \
       p =g_{yz} \quad \text{where} g=2^{12}(x-x^2)^2 (y-y^2)^2 (z-z^2)^2 .
}
We use tetrahedral meshes shown in Figure \ref{grid3}.
The results of the 3D $P_k$-$P_{k+1}$ weak Galerkin finite element methods
    are listed in Table \ref{t3}.
The results show that the method is stable and is of two-order superconvergence (for velocity).

\begin{table}[ht]
  \centering   \renewcommand{\arraystretch}{1.05}
  \caption{Error profiles for solution \eqref{s2} on wedge grids shown in Figure \ref{grid3}. }
\label{t3}
\begin{tabular}{c|cc|cc|cc}
\hline
Grid & $\|Q_h \bu- \bu_h \|_0 $  &rate &  $\3bar Q_h \bu- \bu_h \3bar $ &rate
  &  $\|p - p_h \|_0 $ &rate   \\
\hline
 &\multicolumn{6}{c}{by the $P_0^2$-$P_1^2$-$P_1$ WG finite element } \\ \hline
 4&    0.8167E-01& 1.59&    0.1864E+01& 1.83&    0.5772E+00& 1.78 \\
 5&    0.2228E-01& 1.87&    0.4851E+00& 1.94&    0.1575E+00& 1.87 \\
 6&    0.5689E-02& 1.97&    0.1228E+00& 1.98&    0.3776E-01& 2.06 \\
\hline
 &\multicolumn{6}{c}{by the $P_1^2$-$P_2^2$-$P_2$ WG finite element } \\ \hline
 3&    0.6428E-01& 3.50&    0.4486E+01& 2.52&    0.6305E+00& 3.23 \\
 4&    0.4636E-02& 3.79&    0.6105E+00& 2.88&    0.8163E-01& 2.95 \\
 5&    0.2856E-03& 4.02&    0.7796E-01& 2.97&    0.9492E-02& 3.10 \\
\hline
 &\multicolumn{6}{c}{by the $P_2^2$-$P_3^2$-$P_3$ WG finite element } \\ \hline
 2&    0.7217E+00& 3.28&    0.3793E+02& 1.89&    0.2623E+01& 5.54 \\
 3&    0.2563E-01& 4.82&    0.2898E+01& 3.71&    0.2215E+00& 3.57 \\
 4&    0.8352E-03& 4.94&    0.1942E+00& 3.90&    0.1439E-01& 3.94 \\
 \hline
\end{tabular}%
\end{table}%

\end{document}